\documentclass{amsart}
\usepackage{amssymb}
\newcommand{\pn}{\par\noindent}
\newcommand{\pmn}{\par\medskip\noindent}
\newcommand{\pbn}{\par\bigskip\noindent}
\begin{document}
\title{Combinatorics of generic 5-degree polynomials}
\author{Yury Kochetkov}
\date{}
\begin{abstract} We consider the space $P$ of generic complex 5-degree
polynomials. Critical values of such polynomial, i.e. four points
in the complex plane, either are vertices of a convex quadrangle
$Q$, or vertices of a triangle $T$ with one point inside $T$. The
inverse image of $Q$ is a tree-like connected structure of five
ovals (a cactus). The inverse image of $T$ is also a cactus, but
of four ovals. Transformations of cacti of the first type into
cacti of the second type and vice versa allow one to represent the
space $P$ as a ribbon bipartite graph of genus 3.
\end{abstract}

\email{yukochetkov@hse.ru, yuyukochetkov@gmail.com} \maketitle

\section{Introduction}
\pn Let $U$ be the set of critical values of a complex rational
function $f$ and let $O$ be an oval, such that $U\subset O$. The
inverse image $C=f^{-1}(O)$ is a collection of ovals, where each
two ovals from $C$ are either disjoint, or intersect in critical
points of $f$. W.Thurston (W.Thurston 2010, see also \cite{KL})
gave the description of $C$ in the case of a generic rational $f$
and call $C$ the \emph{shape} of $f$. \pmn Now let $f$ be a
generic complex polynomial $f$ of degree $n$. Here $C$ is a
collection of $n$ ovals, each two of them are either disjoint, or
have exactly one common point --- a critical point of $f$. Thus,
we can consider the following combinatorial object --- a connected
tree-like collection of $n$ ovals, where: a) any two ovals have
either one common point and touch externally, or are disjoint; b)
a common point is a common point of exactly two ovals; c) each
oval have $n-1$ marked labelled points, where the cyclic order of
labels is the same in all ovals; d) a common point of two ovals is
a marked point; e) two different common points have different
labels. Such construction will be called a combinatorial
$n$-cactus. \pmn {\bf Example.} There are two combinatorial
4-cacti:
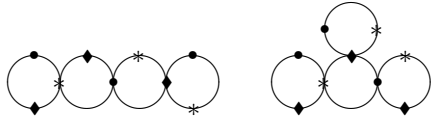
\[\begin{picture}(160,50) \multiput(10,15)(20,0){4}{\circle{20}}
\put(17,12){$\ast$} \put(8,3){\tiny $\blacklozenge$}
\put(10,25){\circle*{3}} \put(40,15){\circle*{3}}
\put(28,23){\tiny $\blacklozenge$} \put(58,13){\tiny
$\blacklozenge$} \put(47,22){$\ast$} \put(68,2){$\ast$}
\put(70,25){\circle*{3}}

\multiput(110,15)(20,0){3}{\circle{20}} \put(117,12){$\ast$}
\put(108,3){\tiny $\blacklozenge$} \put(110,25){\circle*{3}}
\put(140,15){\circle*{3}} \put(128,23){\tiny $\blacklozenge$}
\put(148,3){\tiny $\blacklozenge$} \put(148,22){$\ast$}
\put(130,35){\circle{20}} \put(137,32){$\ast$}
\put(120,35){\circle*{3}}
\end{picture}\] \pmn If $p\in\mathbb{C}[z]$ is a generic polynomial of
degree $n$, then there are $n-1$ points in its set $U$. If $O$ is
a closed Jordan curve that pass through every points of $U$, then
$p^{-1}(O)$ is a combinatorial $n$-cactus. However, different
choices of $O$ can give us all combinatorial cacti. Thus, we must
study those situations, where the curve $O$ can be defined in some
unique way.  \pmn \emph{Remark.} An obvious solution is to study
the space M of degree $n$ complex polynomials with exactly three
critical values. Then we can represent $M$ as a 3-valent graph
(see \cite{Ko}). But such polynomials are not generic.\pmn Let now
$p\in\mathbb{C}[z]$ be a generic 5-degree polynomial. Its critical
points $z_1,z_2,z_3,z_4$ are pairwise different and have
multiplicity 2. Critical values $u_1,u_2,u_3,u_4$, $u_i=p(z_i)$,
$i=1,\ldots,4$, are pairwise different and no three of them belong
to one line. Either points $u_i$, $i=1,\ldots,4$, are vertices of
a convex quadrangle $Q$ --- the first case, or three points are
vertices of a triangle $T$ and the forth is inside it --- the
second case. \subsection{The first case} Here the inverse image
$p^{-1}(Q)$ is a combinatorial 5-cactus. \pmn We will denote
points $u_1,u_2,u_3,u_4$ (and their inverse images) by symbols
$\ast,\circ,\bullet$ and $\diamond$, respectively, and we will fix
the cyclic order of points $u_i$ in $Q$:
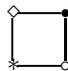
\[\begin{picture}(20,30) \multiput(2,5)(0,20){2}{\line(1,0){16}}
\multiput(0,7)(20,0){2}{\line(0,1){16}} \put(-2,3){$\ast$}
\put(20,5){\circle{3}} \put(20,25){\circle*{3}}
\put(-2,23){$\diamond$} \end{picture}\]  There are eight different cacti
of the first type:
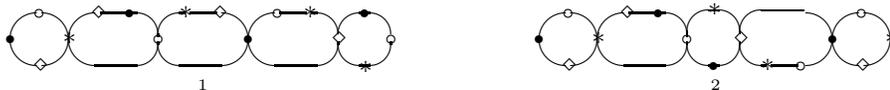
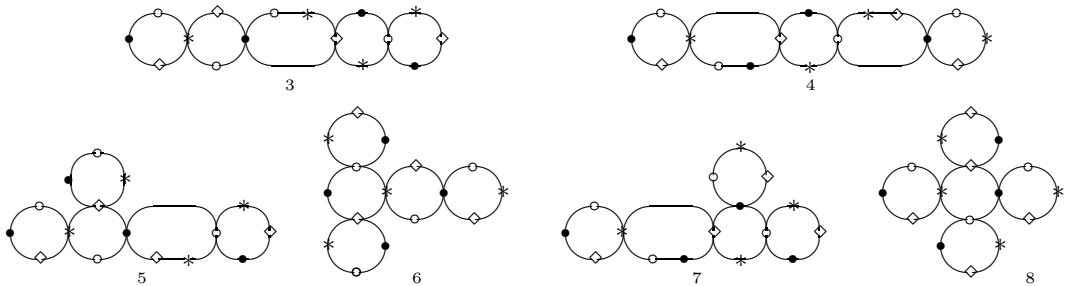
\[\begin{picture}(335,35) \put(10,20){\oval(20,20)[l]}
\put(12,20){\oval(20,20)[r]} \put(0,20){\circle*{3}}
\put(20,18){$\ast$} \put(11,30){\circle{3}} \put(9,8){$\diamond$}

\put(32,20){\oval(20,20)[l]} \put(31,28){$\diamond$}
\put(34,30){\line(1,0){10}} \put(45,30){\circle*{3}}
\put(32,10){\line(1,0){14}} \put(46,19){\oval(20,18)[rb]}
\put(46,21){\oval(20,18)[rt]} \put(56,20){\circle{3}}

\put(66,21){\oval(20,18)[lt]} \put(66,19){\oval(20,18)[lb]}
\put(64,28){$\ast$} \put(68,30){\line(1,0){10}}
\put(66,10){\line(1,0){14}} \put(77,28){$\diamond$}
\put(80,20){\oval(20,20)[r]} \put(90,20){\circle*{3}}

\put(100,20){\oval(20,20)[l]} \put(101,30){\circle{3}}
\put(102,30){\line(1,0){10}} \put(111,28){$\ast$}
\put(100,10){\line(1,0){14}} \put(114,22){\oval(20,16)[rt]}
\put(114,19){\oval(20,18)[rb]} \put(122,18){$\diamond$}

\put(134,22){\oval(20,16)[t]} \put(134,19){\oval(20,18)[b]}
\put(144,20){\circle{3}} \put(134,30){\circle*{3}}
\put(132,7){$\ast$}

\put(71,1){\tiny 1}

\put(210,20){\oval(20,20)[l]} \put(212,20){\oval(20,20)[r]}
\put(200,20){\circle*{3}} \put(220,18){$\ast$}
\put(211,30){\circle{3}} \put(209,8){$\diamond$}

\put(232,20){\oval(20,20)[l]} \put(231,28){$\diamond$}
\put(234,30){\line(1,0){10}} \put(245,30){\circle*{3}}
\put(232,10){\line(1,0){14}} \put(246,19){\oval(20,18)[rb]}
\put(246,21){\oval(20,18)[rt]} \put(256,20){\circle{3}}

\put(266,22){\oval(20,18)[t]} \put(266,19){\oval(20,18)[b]}
\put(274,18){$\diamond$} \put(264,29){$\ast$}
\put(266,10){\circle*{3}}

\put(286,22){\oval(20,18)[lt]} \put(286,19){\oval(20,18)[lb]}
\put(284,8){$\ast$} \put(288,10){\line(1,0){10}}
\put(299,10){\circle{3}} \put(286,31){\line(1,0){14}}
\put(301,20){\oval(20,20)[r]} \put(311,20){\circle*{3}}

\put(321,20){\oval(20,20)[l]} \put(320,8){$\diamond$}
\put(322,30){\circle{3}} \put(323,20){\oval(20,20)[r]}
\put(331,18){$\ast$}

\put(265,1){\tiny 2}
\end{picture}\]

\[\begin{picture}(310,35) \put(10,20){\oval(20,20)[l]}
\put(12,20){\oval(20,20)[r]} \put(0,20){\circle*{3}}
\put(20,18){$\ast$} \put(11,30){\circle{3}} \put(9,8){$\diamond$}

\put(32,20){\oval(20,20)[l]} \put(33,10){\circle{3}}
\put(31,28){$\diamond$} \put(34,20){\oval(20,20)[r]}
\put(44,20){\circle*{3}}

\put(54,20){\oval(20,20)[l]} \put(55,30){\circle{3}}
\put(56,30){\line(1,0){10}} \put(65,27){$\ast$}
\put(54,10){\line(1,0){14}} \put(68,21){\oval(20,18)[rt]}
\put(68,19){\oval(20,18)[rb]} \put(76,18){$\diamond$}

\put(88,21){\oval(20,18)[t]} \put(88,19){\oval(20,18)[b]}
\put(98,20){\circle{3}} \put(88,30){\circle*{3}}
\put(86,8){$\ast$}

\put(108,21){\oval(20,18)[t]} \put(108,19){\oval(20,18)[b]}
\put(116,18){$\diamond$} \put(108,10){\circle*{3}}
\put(106,28){$\ast$}

\put(59,1){\tiny 3}

\put(200,20){\oval(20,20)[l]} \put(202,20){\oval(20,20)[r]}
\put(190,20){\circle*{3}} \put(210,18){$\ast$}
\put(201,30){\circle{3}} \put(199,8){$\diamond$}

\put(222,20){\oval(20,20)[l]} \put(223,10){\circle{3}}
\put(224,10){\line(1,0){10}} \put(235,10){\circle*{3}}
\put(222,30){\line(1,0){14}} \put(236,21){\oval(20,18)[rt]}
\put(236,19){\oval(20,18)[rb]} \put(244,18){$\diamond$}

\put(256,21){\oval(20,18)[lt]} \put(256,19){\oval(20,18)[lb]}
\put(255,7){$\ast$} \put(257,30){\circle*{3}}
\put(258,21){\oval(20,18)[rt]} \put(258,19){\oval(20,18)[rb]}
\put(268,20){\circle{3}}

\put(278,21){\oval(20,18)[lt]} \put(278,19){\oval(20,18)[lb]}
\put(277,27){$\ast$} \put(280,30){\line(1,0){10}}
\put(288,27){$\diamond$} \put(278,10){\line(1,0){14}}
\put(292,20){\oval(20,20)[r]} \put(302,20){\circle*{3}}

\put(312,20){\oval(20,20)[l]} \put(313,30){\circle{3}}
\put(311,8){$\diamond$} \put(314,20){\oval(20,20)[r]}
\put(322,18){$\ast$} \put(256,1){\tiny 4} \end{picture}\]

\[\begin{picture}(400,50) \put(10,20){\oval(20,20)[l]}
\put(12,20){\oval(20,20)[r]} \put(0,20){\circle*{3}}
\put(20,18){$\ast$} \put(11,30){\circle{3}} \put(9,8){$\diamond$}

\put(32,20){\oval(20,20)[l]} \put(33,10){\circle{3}}
\put(31,28){$\diamond$} \put(34,20){\oval(20,20)[r]}
\put(44,20){\circle*{3}}

\put(32,40){\oval(18,20)[l]} \put(34,40){\oval(18,20)[r]}
\put(33,50){\circle{3}} \put(22,40){\circle*{3}}
\put(41,38){$\ast$}

\put(54,20){\oval(20,20)[l]} \put(53,8){$\diamond$}
\put(56,10){\line(1,0){10}} \put(65,7){$\ast$}
\put(54,30){\line(1,0){14}} \put(68,21){\oval(20,18)[rt]}
\put(68,19){\oval(20,18)[rb]} \put(78,20){\circle{3}}

\put(88,21){\oval(20,18)[t]} \put(88,19){\oval(20,18)[b]}
\put(96,18){$\diamond$} \put(88,10){\circle*{3}}
\put(86,28){$\ast$}

\put(48,1){\tiny 5}

\put(130,35){\oval(20,20)[l]} \put(131,45){\circle{3}}
\put(129,23){$\diamond$} \put(132,35){\oval(20,20)[r]}
\put(120,35){\circle*{3}} \put(140,33){$\ast$}

\put(130,15){\oval(20,20)[l]} \put(118,13){$\ast$}
\put(131,5){\circle{3}} \put(132,15){\oval(20,20)[r]}
\put(131,5){\circle{3}} \put(142,15){\circle*{3}}

\put(130,55){\oval(20,20)[l]} \put(129,63){$\diamond$}
\put(118,53){$\ast$} \put(132,55){\oval(20,20)[r]}
\put(142,55){\circle*{3}}

\put(152,35){\oval(20,20)[l]} \put(153,25){\circle{3}}
\put(151,43){$\diamond$} \put(154,35){\oval(20,20)[r]}
\put(164,35){\circle*{3}}

\put(174,35){\oval(20,20)[l]} \put(173,23){$\diamond$}
\put(175,45){\circle{3}} \put(176,35){\oval(20,20)[r]}
\put(184,33){$\ast$}

\put(152,1){\tiny 6}

\put(220,20){\oval(20,20)[l]} \put(222,20){\oval(20,20)[r]}
\put(210,20){\circle*{3}} \put(229,18){$\ast$}
\put(221,30){\circle{3}} \put(219,8){$\diamond$}

\put(242,20){\oval(20,20)[l]} \put(243,10){\circle{3}}
\put(244,10){\line(1,0){10}} \put(255,10){\circle*{3}}
\put(242,30){\line(1,0){14}} \put(256,21){\oval(20,18)[rt]}
\put(256,19){\oval(20,18)[rb]} \put(264,18){$\diamond$}

\put(276,21){\oval(20,18)[t]} \put(276,19){\oval(20,18)[b]}
\put(286,20){\circle{3}} \put(276,30){\circle*{3}}
\put(274,7){$\ast$}

\put(276,40){\oval(20,20)[b]} \put(276,42){\oval(20,20)[t]}
\put(266,41){\circle{3}} \put(284,39){$\diamond$}
\put(274,50){$\ast$}

\put(296,21){\oval(20,18)[t]} \put(296,19){\oval(20,18)[b]}
\put(304,18){$\diamond$} \put(296,10){\circle*{3}}
\put(294,28){$\ast$}

\put(258,1){\tiny 7}

\put(340,35){\oval(20,20)[l]} \put(341,45){\circle{3}}
\put(339,23){$\diamond$} \put(342,35){\oval(20,20)[r]}
\put(330,35){\circle*{3}} \put(350,33){$\ast$}

\put(362,35){\oval(20,20)[l]} \put(364,35){\oval(20,20)[r]}
\put(363,25){\circle{3}} \put(361,43){$\diamond$}
\put(374,35){\circle*{3}}

\put(362,15){\oval(20,20)[l]} \put(364,15){\oval(20,20)[r]}
\put(361,3){$\diamond$} \put(352,15){\circle*{3}}
\put(372,13){$\ast$}

\put(362,55){\oval(20,20)[l]} \put(364,55){\oval(20,20)[r]}
\put(361,63){$\diamond$} \put(374,55){\circle*{3}}
\put(350,53){$\ast$}

\put(384,35){\oval(20,20)[l]} \put(386,35){\oval(20,20)[r]}
\put(383,23){$\diamond$} \put(385,45){\circle{3}}
\put(394,33){$\ast$}

\put(384,1){\tiny 8}
\end{picture}\]

\subsection{The second case} Here points $u_1,u_2,u_3$ are
vertices of the triangle $T$ and the point $u_4$ is inside $T$.
The inverse image $p^{-1}(T)$ is a cactus of four ovals --- three
"small" and one "big". The border of each small oval contains an
inverse image of $u_1$, an inverse image of $u_2$ and an inverse
image of $u_3$, and inside each small oval is an inverse image of
$u_4$ (not $z_4$). The border of the big oval contains six points:
two inverse images of $u_1$, two inverse images of $u_2$ and two
inverse images of $u_3$. Point $z_4$ is inside the big oval. There
are nine cacti of the second type:
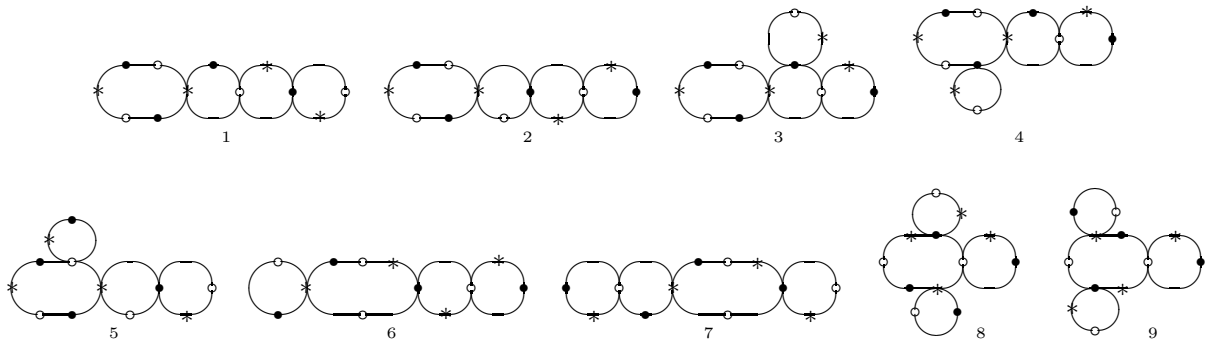
\[\begin{picture}(385,50) \put(10,20){\oval(20,20)[l]}
\put(-2,18){$\ast$} \put(11,30){\circle*{3}}
\put(11,10){\circle{3}} \multiput(12,10)(0,20){2}{\line(1,0){10}}
\put(23,10){\circle*{3}} \put(23,30){\circle{3}}
\put(24,20){\oval(20,20)[r]} \put(32,18){$\ast$}

\put(44,21){\oval(20,18)[t]} \put(44,19){\oval(20,18)[b]}
\put(54,20){\circle{3}} \put(44,30){\circle*{3}}

\put(64,21){\oval(20,18)[t]} \put(64,19){\oval(20,18)[b]}
\put(74,20){\circle*{3}} \put(62,27){$\ast$}

\put(84,21){\oval(20,18)[t]} \put(84,19){\oval(20,18)[b]}
\put(94,20){\circle{3}} \put(82,8){$\ast$}

\put(47,1){\tiny 1}

\put(120,20){\oval(20,20)[l]} \put(108,18){$\ast$}
\put(121,30){\circle*{3}} \put(121,10){\circle{3}}
\multiput(122,10)(0,20){2}{\line(1,0){10}}
\put(133,10){\circle*{3}} \put(133,30){\circle{3}}
\put(134,20){\oval(20,20)[r]} \put(142,18){$\ast$}

\put(154,20){\oval(20,20)[t]} \put(153,20){\oval(18,20)[lb]}
\put(155,20){\oval(18,20)[rb]} \put(154,10){\circle{3}}
\put(164,20){\circle*{3}}

\put(174,21){\oval(20,18)[t]} \put(174,19){\oval(20,18)[b]}
\put(184,20){\circle{3}} \put(172,7){$\ast$}

\put(194,21){\oval(20,18)[t]} \put(194,19){\oval(20,18)[b]}
\put(204,20){\circle*{3}} \put(192,27){$\ast$}

\put(161,1){\tiny 2}

\put(230,20){\oval(20,20)[l]} \put(218,18){$\ast$}
\put(231,30){\circle*{3}} \put(231,10){\circle{3}}
\multiput(232,10)(0,20){2}{\line(1,0){10}}
\put(243,10){\circle*{3}} \put(243,30){\circle{3}}
\put(244,20){\oval(20,20)[r]} \put(252,18){$\ast$}

\put(264,21){\oval(20,18)[t]} \put(264,19){\oval(20,18)[b]}
\put(274,20){\circle{3}} \put(264,30){\circle*{3}}

\put(284,21){\oval(20,18)[t]} \put(284,19){\oval(20,18)[b]}
\put(294,20){\circle*{3}} \put(282,27){$\ast$}

\put(263,40){\oval(18,20)[l]} \put(265,40){\oval(18,20)[r]}
\put(264,50){\circle{3}} \put(272,38){$\ast$}

\put(256,1){\tiny 3}

\put(320,40){\oval(20,20)[l]} \put(308,38){$\ast$}
\put(321,50){\circle*{3}} \put(321,30){\circle{3}}
\multiput(322,30)(0,20){2}{\line(1,0){10}}
\put(333,30){\circle*{3}} \put(333,50){\circle{3}}
\put(334,40){\oval(20,20)[r]} \put(342,38){$\ast$}

\put(354,41){\oval(20,18)[t]} \put(354,39){\oval(20,18)[b]}
\put(364,40){\circle{3}} \put(354,50){\circle*{3}}

\put(374,41){\oval(20,18)[t]} \put(374,39){\oval(20,18)[b]}
\put(384,40){\circle*{3}} \put(372,48){$\ast$}

\put(332,21){\oval(16,16)[l]} \put(334,21){\oval(16,16)[r]}
\put(333,13){\circle{3}} \put(322,18){$\ast$}

\put(347,1){\tiny 4}
\end{picture}\]

\[\begin{picture}(450,45)
\put(10,20){\oval(20,20)[l]} \put(-2,18){$\ast$}
\put(11,30){\circle*{3}} \put(11,10){\circle{3}}
\multiput(12,10)(0,20){2}{\line(1,0){10}} \put(23,10){\circle*{3}}
\put(23,30){\circle{3}} \put(24,20){\oval(20,20)[r]}
\put(32,18){$\ast$}

\put(44,20){\oval(20,20)[lb]} \put(46,20){\oval(20,20)[rb]}
\put(45,20){\oval(22,20)[t]} \put(45,10){\circle{3}}
\put(56,20){\circle*{3}}

\put(66,21){\oval(20,18)[t]} \put(66,19){\oval(20,18)[b]}
\put(76,20){\circle{3}} \put(64,7){$\ast$}

\put(22,38){\oval(16,16)[l]} \put(24,38){\oval(16,16)[r]}
\put(23,46){\circle*{3}} \put(12,36){$\ast$}

\put(37,1){\tiny 5}

\put(100,20){\oval(20,20)[l]} \put(102,20){\oval(20,20)[r]}
\put(101,30){\circle{3}} \put(101,10){\circle*{3}}
\put(109,18){$\ast$}

\put(122,20){\oval(20,20)[l]} \put(122,30){\circle*{3}}
\multiput(122,10)(0,20){2}{\line(1,0){10}}
\multiput(133,10)(0,20){2}{\circle{3}}
\multiput(134,10)(0,20){2}{\line(1,0){10}}
\put(144,20){\oval(20,20)[r]} \put(142,27){$\ast$}
\put(154,20){\circle*{3}}

\put(164,21){\oval(20,18)[t]} \put(164,19){\oval(20,18)[b]}
\put(174,20){\circle{3}} \put(162,8){$\ast$}

\put(184,21){\oval(20,18)[t]} \put(184,19){\oval(20,18)[b]}
\put(194,20){\circle*{3}} \put(182,28){$\ast$}

\put(142,1){\tiny 6}

\put(220,21){\oval(20,18)[t]} \put(220,19){\oval(20,18)[b]}
\put(210,20){\circle*{3}} \put(230,20){\circle{3}}
\put(218,7){$\ast$}

\put(240,21){\oval(20,18)[t]} \put(240,19){\oval(20,18)[b]}
\put(248,18){$\ast$} \put(240,10){\circle*{3}}

\put(260,20){\oval(20,20)[l]} \put(260,30){\circle*{3}}
\multiput(260,10)(0,20){2}{\line(1,0){10}}
\multiput(271,10)(0,20){2}{\circle{3}}
\multiput(272,10)(0,20){2}{\line(1,0){10}}
\put(282,20){\oval(20,20)[r]} \put(292,20){\circle*{3}}
\put(280,27){$\ast$}

\put(302,21){\oval(20,18)[t]} \put(302,19){\oval(20,18)[b]}
\put(312,20){\circle{3}} \put(300,7){$\ast$}

\put(262,1){\tiny 7}

\put(340,31){\oval(20,18)[lt]} \put(340,29){\oval(20,18)[lb]}
\put(330,30){\circle{3}}
\multiput(340,20)(0,20){2}{\line(1,0){10}}
\put(350,31){\oval(20,18)[rt]} \put(350,29){\oval(20,18)[rb]}
\put(360,30){\circle{3}} \multiput(340,20)(10,20){2}{\circle*{3}}
\multiput(338,37)(10,-20){2}{$\ast$}

\put(370,31){\oval(20,18)[t]} \put(370,29){\oval(20,18)[b]}
\put(380,30){\circle*{3}} \put(368,37){$\ast$}

\put(349,48){\oval(16,16)[l]} \put(351,48){\oval(16,16)[r]}
\put(350,56){\circle{3}} \put(357,46){$\ast$}

\put(350,12){\oval(16,16)[t]} \put(350,10){\oval(16,16)[b]}
\put(342,11){\circle{3}} \put(358,11){\circle*{3}}

\put(365,1){\tiny 8}

\put(410,31){\oval(20,18)[lt]} \put(410,29){\oval(20,18)[lb]}
\put(400,30){\circle{3}}
\multiput(410,20)(0,20){2}{\line(1,0){10}}
\put(420,31){\oval(20,18)[rt]} \put(420,29){\oval(20,18)[rb]}
\put(430,30){\circle{3}} \multiput(410,20)(10,20){2}{\circle*{3}}
\multiput(408,37)(10,-20){2}{$\ast$}

\put(440,31){\oval(20,18)[t]} \put(440,29){\oval(20,18)[b]}
\put(450,30){\circle*{3}} \put(438,37){$\ast$}

\put(409,12){\oval(16,16)[l]} \put(411,12){\oval(16,16)[r]}
\put(410,4){\circle{3}} \put(399,10){$\ast$}

\put(410,48){\oval(16,16)[b]} \put(410,50){\oval(16,16)[t]}
\put(402,49){\circle*{3}} \put(418,49){\circle{3}}

\put(430,1){\tiny 9}
\end{picture}\]  Now we will describe how to transform a cactus of the
first type into a cactus of the second type and vice versa.

\section{Transformation 1}
\pn Each first type cactus in four ways can be transformed into a
second type cactus. We choose a vertex in the quadrangle $Q$ and
connect the previous and the next vertices by an arc in the
exterior of $Q$:
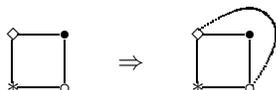
\[\begin{picture}(100,40) \multiput(2,5)(0,20){2}{\line(1,0){16}}
\multiput(0,7)(20,0){2}{\line(0,1){16}} \put(-2,3){$\ast$}
\put(20,5){\circle{3}} \put(20,25){\circle*{3}}
\put(-2,23){$\diamond$}

\put(40,12){$\Rightarrow$}

\multiput(72,5)(0,20){2}{\line(1,0){16}}
\multiput(70,7)(20,0){2}{\line(0,1){16}} \put(68,3){$\ast$}
\put(90,5){\circle{3}} \put(90,25){\circle*{3}}
\put(68,23){$\diamond$} \qbezier[35](71,27)(90,40)(97,32)
\qbezier[35](91,6)(105,25)(97,32)
\end{picture}\] Thus obtained figure can be considered as a
"triangle" $\widetilde{T}$, where the chosen vertex is inside it.
The inverse image of $\widetilde{T}$ is a cactus of the second
type. For example,
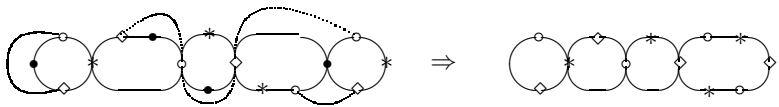
\[\begin{picture}(300,50) \put(20,20){\oval(20,20)[l]}
\put(22,20){\oval(20,20)[r]} \put(10,20){\circle*{3}}
\put(30,18){$\ast$} \put(21,30){\circle{3}} \put(19,8){$\diamond$}

\put(42,20){\oval(20,20)[l]} \put(41,28){$\diamond$}
\put(44,30){\line(1,0){10}} \put(55,30){\circle*{3}}
\put(42,10){\line(1,0){14}} \put(56,19){\oval(20,18)[rb]}
\put(56,21){\oval(20,18)[rt]} \put(66,20){\circle{3}}

\put(76,22){\oval(20,18)[t]} \put(76,19){\oval(20,18)[b]}
\put(84,18){$\diamond$} \put(74,29){$\ast$}
\put(76,10){\circle*{3}}

\put(96,22){\oval(20,18)[lt]} \put(96,19){\oval(20,18)[lb]}
\put(94,8){$\ast$} \put(98,10){\line(1,0){10}}
\put(109,10){\circle{3}} \put(96,31){\line(1,0){14}}
\put(111,20){\oval(20,20)[r]} \put(121,20){\circle*{3}}

\put(131,20){\oval(20,20)[l]} \put(130,8){$\diamond$}
\put(132,30){\circle{3}} \put(133,20){\oval(20,20)[r]}
\put(141,18){$\ast$}

\qbezier[40](0,20)(0,35)(19,31) \qbezier[40](0,20)(0,5)(19,10)
\qbezier[40](43,31)(66,50)(66,22)
\qbezier[50](86,22)(86,55)(131,31)
\qbezier[30](131,9)(120,0)(110,9) \qbezier[20](66,18)(66,5)(76,5)
\qbezier[20](76,5)(86,5)(86,18)

\put(160,18){$\Rightarrow$}

\put(200,20){\oval(20,20)[l]}\put(201,30){\circle{3}}
\put(199,8){$\diamond$} \put(202,20){\oval(20,20)[r]}
\put(210,18){$\ast$}

\put(223,19){\oval(22,18)[b]} \put(222,21){\oval(20,18)[lt]}
\put(221,27){$\diamond$} \put(224,21){\oval(20,18)[rt]}
\put(234,20){\circle{3}}

\put(244,21){\oval(20,18)[t]} \put(244,19){\oval(20,18)[b]}
\put(252,18){$\diamond$} \put(241,27){$\ast$}

\put(264,21){\oval(20,18)[lt]} \put(264,19){\oval(20,18)[lb]}
\put(265,30){\circle{3}} \put(263,7){$\ast$}
\multiput(266,10)(0,20){2}{\line(1,0){10}}
\put(277,10){\circle{3}} \put(275,27){$\ast$}
\put(278,21){\oval(20,18)[rt]} \put(278,19){\oval(20,18)[rb]}
\put(286,18){$\diamond$}
\end{picture}\] Thus, given a first type cactus, we can construct four
second type cacti. Moreover, a cyclic order of vertices in the
given first type cactus induces the same cyclic order in the set
of (constructed from it) second type cacti. In the following
figures the number of a first type cactus is in a circle and
numbers of constructed cacti --- in squares (in the counter
clockwise order).
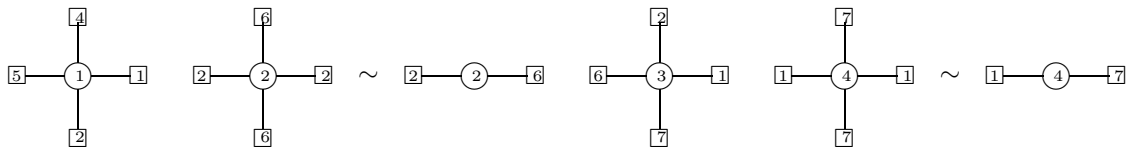
\[\begin{picture}(420,55) \put(-2,27){$\square$}
\put(5,30){\line(1,0){15}} \put(25,30){\circle{10}}
\put(30,30){\line(1,0){15}} \put(25,35){\line(0,1){15}}
\put(25,25){\line(0,-1){15}} \put(21,49){$\square$}
\put(44,27){$\square$} \put(21,3){$\square$}

\put(0,28){\tiny 5} \put(24,50){\tiny 4} \put(47,28){\tiny 1}
\put(24,5){\tiny 2} \put(24,28){\tiny 1}

\put(68,27){$\square$} \put(75,30){\line(1,0){15}}
\put(95,30){\circle{10}} \put(100,30){\line(1,0){15}}
\put(95,35){\line(0,1){15}} \put(95,25){\line(0,-1){15}}
\put(91,49){$\square$} \put(114,27){$\square$}
\put(91,3){$\square$}

\put(70,28){\tiny 2} \put(94,50){\tiny 6} \put(117,28){\tiny 2}
\put(94,5){\tiny 6} \put(94,28){\tiny 2}

\put(131,28){$\sim$}

\put(148,27){$\square$} \put(155,30){\line(1,0){15}}
\put(175,30){\circle{10}} \put(180,30){\line(1,0){15}}
\put(194,27){$\square$} \put(150,28){\tiny 2} \put(197,28){\tiny
6} \put(174,28){\tiny 2}

\put(218,27){$\square$} \put(225,30){\line(1,0){15}}
\put(245,30){\circle{10}} \put(250,30){\line(1,0){15}}
\put(245,35){\line(0,1){15}} \put(245,25){\line(0,-1){15}}
\put(241,49){$\square$} \put(264,27){$\square$}
\put(241,3){$\square$}

\put(220,28){\tiny 6} \put(244,50){\tiny 2} \put(267,28){\tiny 1}
\put(244,5){\tiny 7} \put(244,28){\tiny 3}

\put(288,27){$\square$} \put(295,30){\line(1,0){15}}
\put(315,30){\circle{10}} \put(320,30){\line(1,0){15}}
\put(315,35){\line(0,1){15}} \put(315,25){\line(0,-1){15}}
\put(311,49){$\square$} \put(334,27){$\square$}
\put(311,3){$\square$}

\put(290,28){\tiny 1} \put(314,50){\tiny 7} \put(337,28){\tiny 1}
\put(314,5){\tiny 7} \put(314,28){\tiny 4}

\put(351,28){$\sim$}

\put(368,27){$\square$} \put(375,30){\line(1,0){15}}
\put(395,30){\circle{10}} \put(400,30){\line(1,0){15}}
\put(414,27){$\square$} \put(370,28){\tiny 1} \put(417,28){\tiny
7} \put(394,28){\tiny 4}
\end{picture}\]

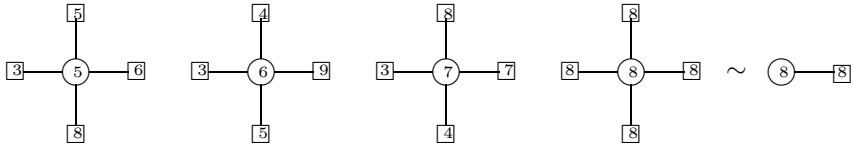
\[\begin{picture}(420,55) \put(-2,27){$\square$}
\put(5,30){\line(1,0){15}} \put(25,30){\circle{10}}
\put(30,30){\line(1,0){15}} \put(25,35){\line(0,1){15}}
\put(25,25){\line(0,-1){15}} \put(21,49){$\square$}
\put(44,27){$\square$} \put(21,3){$\square$}

\put(1,29){\tiny 3} \put(24,50){\tiny 5} \put(47,29){\tiny 6}
\put(24,5){\tiny 8} \put(24,28){\tiny 5}

\put(68,27){$\square$} \put(75,30){\line(1,0){15}}
\put(95,30){\circle{10}} \put(100,30){\line(1,0){15}}
\put(95,35){\line(0,1){15}} \put(95,25){\line(0,-1){15}}
\put(91,49){$\square$} \put(114,27){$\square$}
\put(91,3){$\square$}

\put(71,29){\tiny 3} \put(94,50){\tiny 4} \put(117,29){\tiny 9}
\put(94,5){\tiny 5} \put(94,28){\tiny 6}

\put(138,27){$\square$} \put(145,30){\line(1,0){15}}
\put(165,30){\circle{10}} \put(170,30){\line(1,0){15}}
\put(165,35){\line(0,1){15}} \put(165,25){\line(0,-1){15}}
\put(161,49){$\square$} \put(184,27){$\square$}
\put(161,3){$\square$}

\put(140,29){\tiny 3} \put(164,50){\tiny 8} \put(187,29){\tiny 7}
\put(164,5){\tiny 4} \put(164,28){\tiny 7}

\put(208,27){$\square$} \put(215,30){\line(1,0){15}}
\put(235,30){\circle{10}} \put(240,30){\line(1,0){15}}
\put(235,35){\line(0,1){15}} \put(235,25){\line(0,-1){15}}
\put(231,49){$\square$} \put(254,27){$\square$}
\put(231,3){$\square$}

\put(210,29){\tiny 8} \put(234,50){\tiny 8} \put(257,29){\tiny 8}
\put(234,5){\tiny 8} \put(234,28){\tiny 8}

\put(271,28){$\sim$}

\put(292,30){\circle{10}} \put(297,30){\line(1,0){15}}
\put(311,26){$\square$} \put(291,28){\tiny 8} \put(313,28){\tiny
8} \end{picture}\]

\section{Transformation 2}
\pn Each second type cactus in three ways can be transformed into
first type cacti: we delete an edge in the triangle $T$ and
endpoints of this edge connect by arcs with the point $u_i$
inside:
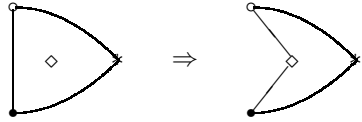
\[\begin{picture}(130,50) \put(0,5){\circle*{3}}
\put(0,45){\circle{3}} \put(0,5){\line(0,1){39}}
\qbezier(0,5)(20,5)(40,25) \qbezier(1,45)(20,45)(40,25)
\put(37,23){$\ast$} \put(12,22){$\diamond$}

\put(60,23){$\Rightarrow$}

\put(90,5){\circle*{3}} \put(90,45){\circle{3}}
\qbezier(90,5)(110,5)(130,25) \qbezier(91,45)(110,45)(130,25)
\put(127,23){$\ast$} \put(103,22){$\diamond$}
\put(90,5){\line(4,5){14}} \put(91,44){\line(4,-5){14}}
\end{picture}\] The corresponding transformation of the cactus is performed
in the following way: in each oval we delete an arc between
vertices $\circ$ and $\bullet$ and these vertices we connect with
the vertex $\diamond$ inside:
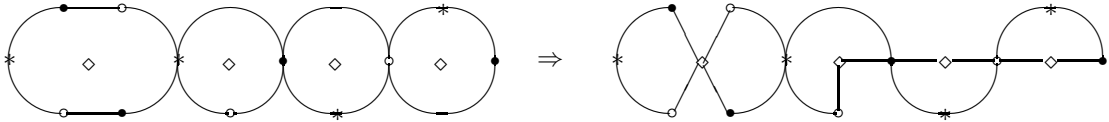
\[\begin{picture}(415,50) \put(20,25){\oval(40,40)[l]}
\put(-2,23){$\ast$} \put(21,45){\circle*{3}}
\put(21,5){\circle{3}} \multiput(22,5)(0,40){2}{\line(1,0){20}}
\put(43,45){\circle{3}} \put(43,5){\circle*{3}}
\put(44,25){\oval(40,40)[r]} \put(62,23){$\ast$}
\put(28,21){$\diamond$}

\put(84,25){\oval(40,40)[t]} \put(83,25){\oval(38,40)[lb]}
\put(85,25){\oval(38,40)[rb]} \put(84,5){\circle{3}}
\put(104,25){\circle*{3}} \put(81,21){$\diamond$}

\put(124,26){\oval(40,38)[t]} \put(124,24){\oval(40,38)[b]}
\put(144,25){\circle{3}} \put(122,2){$\ast$}
\put(121,21){$\diamond$}

\put(164,26){\oval(40,38)[t]} \put(164,24){\oval(40,38)[b]}
\put(184,25){\circle*{3}} \put(162,42){$\ast$}
\put(161,21){$\diamond$}

\put(200,23){$\Rightarrow$}

\put(250,25){\oval(40,40)[l]} \put(228,23){$\ast$}
\put(251,45){\circle*{3}} \put(251,5){\circle{3}}
\put(273,45){\circle{3}} \put(273,5){\circle*{3}}
\put(274,25){\oval(40,40)[r]} \put(292,23){$\ast$}
\put(260,22){$\diamond$}

\put(314,25){\oval(40,40)[t]} \put(313,25){\oval(38,40)[lb]}
\put(314,5){\circle{3}} \put(334,25){\circle*{3}}
\put(312,22){$\diamond$}

\put(354,24){\oval(40,38)[b]} \put(374,25){\circle{3}}
\put(352,2){$\ast$} \put(352,22){$\diamond$}

\put(394,26){\oval(40,38)[t]} \put(414,25){\circle*{3}}
\put(392,42){$\ast$} \put(392,22){$\diamond$}

\put(251,45){\line(1,-2){10}} \put(252,6){\line(1,2){9}}
\put(273,44){\line(-1,-2){10}} \put(272,5){\line(-1,2){9}}

\put(314,6){\line(0,1){17}} \put(315,25){\line(1,0){18}}

\put(334,25){\line(1,0){17}} \put(373,25){\line(-1,0){16}}
\put(375,25){\line(1,0){16}} \put(414,25){\line(-1,0){17}}
\end{picture}\] Thus, given a second type cactus, we can construct
three first type cacti. Moreover, a cyclic order of vertices in
the given second type cactus induces the same cyclic order in the
set of (constructed from it) first type cacti. In the following
figures the number of a second type cactus is in a square and
numbers of constructed cacti --- in circles (in the counter
clockwise order).
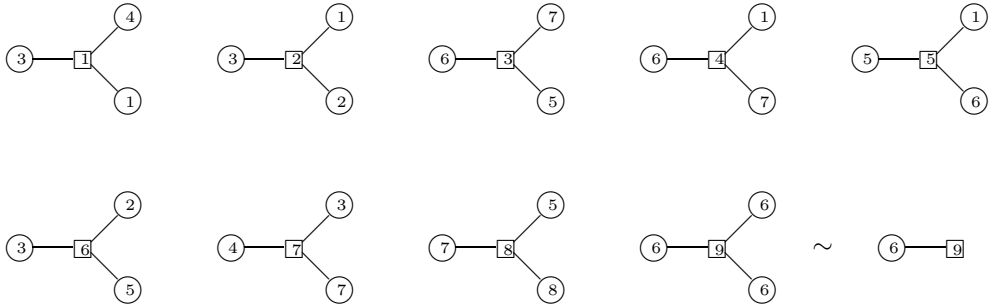
\[\begin{picture}(365,50) \put(0,25){\circle{10}}
\put(5,25){\line(1,0){15}} \put(20,21){$\square$}
\put(27,27){\line(1,1){10}} \put(27,23){\line(1,-1){10}}
\put(41,41){\circle{10}} \put(41,9){\circle{10}}

\put(-1,23){\tiny 3} \put(23,23){\tiny 1} \put(40,39){\tiny 4}
\put(40,7){\tiny 1}

\put(80,25){\circle{10}} \put(85,25){\line(1,0){15}}
\put(100,21){$\square$} \put(107,27){\line(1,1){10}}
\put(107,23){\line(1,-1){10}} \put(121,41){\circle{10}}
\put(121,9){\circle{10}}

\put(79,23){\tiny 3} \put(103,23){\tiny 2} \put(120,39){\tiny 1}
\put(120,7){\tiny 2}

\put(160,25){\circle{10}} \put(165,25){\line(1,0){15}}
\put(180,21){$\square$} \put(187,27){\line(1,1){10}}
\put(187,23){\line(1,-1){10}} \put(201,41){\circle{10}}
\put(201,9){\circle{10}}

\put(159,23){\tiny 6} \put(183,23){\tiny 3} \put(200,39){\tiny 7}
\put(200,7){\tiny 5}

\put(240,25){\circle{10}} \put(245,25){\line(1,0){15}}
\put(260,21){$\square$} \put(267,27){\line(1,1){10}}
\put(267,23){\line(1,-1){10}} \put(281,41){\circle{10}}
\put(281,9){\circle{10}}

\put(239,23){\tiny 6} \put(263,23){\tiny 4} \put(280,39){\tiny 1}
\put(280,7){\tiny 7}

\put(320,25){\circle{10}} \put(325,25){\line(1,0){15}}
\put(340,21){$\square$} \put(347,27){\line(1,1){10}}
\put(347,23){\line(1,-1){10}} \put(361,41){\circle{10}}
\put(361,9){\circle{10}}

\put(319,23){\tiny 5} \put(343,23){\tiny 5} \put(360,39){\tiny 1}
\put(360,7){\tiny 6}
\end{picture}\]

\[\begin{picture}(365,50) \put(0,25){\circle{10}}
\put(5,25){\line(1,0){15}} \put(20,21){$\square$}
\put(27,27){\line(1,1){10}} \put(27,23){\line(1,-1){10}}
\put(41,41){\circle{10}} \put(41,9){\circle{10}}

\put(-1,23){\tiny 3} \put(23,22){\tiny 6} \put(40,39){\tiny 2}
\put(40,7){\tiny 5}

\put(80,25){\circle{10}} \put(85,25){\line(1,0){15}}
\put(100,21){$\square$} \put(107,27){\line(1,1){10}}
\put(107,23){\line(1,-1){10}} \put(121,41){\circle{10}}
\put(121,9){\circle{10}}

\put(79,23){\tiny 4} \put(103,22){\tiny 7} \put(120,39){\tiny 3}
\put(120,7){\tiny 7}

\put(160,25){\circle{10}} \put(165,25){\line(1,0){15}}
\put(180,21){$\square$} \put(187,27){\line(1,1){10}}
\put(187,23){\line(1,-1){10}} \put(201,41){\circle{10}}
\put(201,9){\circle{10}}

\put(159,23){\tiny 7} \put(183,22){\tiny 8} \put(200,39){\tiny 5}
\put(200,7){\tiny 8}

\put(240,25){\circle{10}} \put(245,25){\line(1,0){15}}
\put(260,21){$\square$} \put(267,27){\line(1,1){10}}
\put(267,23){\line(1,-1){10}} \put(281,41){\circle{10}}
\put(281,9){\circle{10}}

\put(239,23){\tiny 6} \put(263,22){\tiny 9} \put(280,39){\tiny 6}
\put(280,7){\tiny 6}

\put(300,23){$\sim$}

\put(330,25){\circle{10}} \put(335,25){\line(1,0){15}}
\put(350,21){$\square$}

\put(329,23){\tiny 6} \put(353,22){\tiny 9}
\end{picture}\]

\section{The ribbon graph}
\pn Now we can construct a ribbon bipartite graph, where the first
type cacti are white vertices and the second type cacti --- black.
In each vertex of the graph the cyclic order of outgoing edges is
defined. Each closed loop, where edges are to the left of our
walk, defines a face of the graph. Below is the list of closed
loops, where numbers of the first type cacti are in circles and
numbers of the second type cacti --- in squares:
$$\ldots\to[2]\to(1)\to[1]\to(4)\to[7]\to(7)\to[8]\to(8)\to[8]\to(5)\to[6]\to
(2)\to[2]\to(1)\to\ldots$$
\begin{multline*}
\ldots\to[1]\to(1)\to[4]\to(6)\to[3]\to(5)\to[8]\to(7)\to[3]\to(6)
\to[5]\to(1)\to[2]\to(3)\to\\
\to[6]\to(5)\to[5]\to(6)\to[9]\to(6)\to[4]\to(7)\to[7]\to(3)\to[1]
\to(1)\to\ldots\end{multline*} \vspace{1mm}
$$\ldots\to[4]\to(1)\to[5]\to(5)\to[3]\to(7)\to[4]\to(1)\to\ldots$$
\vspace{1mm}
$$\ldots\to[2]\to(2)\to[6]\to(3)\to[7]\to(4)\to[1]\to(3)\to[2]
\to(2)\to\ldots$$ As there are 17 vertices, 25 edges and 4 faces,
then our ribbon graph can be embedded in a surface of genus 3.

\vspace{1cm}

\begin{thebibliography}{99}
\bibitem{KL} S.Koch and T.Lei, \emph{On balanced planar graphs, following
W. Thurston}, Arxiv: 1502.04760.
\bibitem{Ko} Yu. Kochetkov, \emph{Zolotarev polynomials of degree 5, 6 and
7 with simple critical points and their moduli spaces}, Arxiv:
2208.02969.
\end{thebibliography}
\end{document}